\documentclass[12pt]{article}
\usepackage{times}

\usepackage{latexsym}
\usepackage{amsmath}
\usepackage[left=1in,top=1in,right=1in]{geometry}
\usepackage{amssymb,amsmath}
\usepackage{amsthm}
\usepackage{amsmath}
\usepackage{amsfonts}
\newtheorem{dfn}{Definition}[section]

\newtheorem{prf}{Proof}[section]
\newtheorem{cor}{Corollary}[section]
\newtheorem{teo}{Theorem}[section]

\theoremstyle{definition}

\numberwithin{equation}{section}

\usepackage{amsmath,amssymb}
\usepackage{titlesec}
\titlelabel{\thetitle.\quad}
\usepackage{amsthm}

\newcommand{\NCPM}{normal metric contact pair manifold}
\title{\textbf{Some Flatness Conditions on Normal Metric Contact Pairs} }
\author{\.Inan \"Unal$ ^1 $\\ Department of Computer Engineering, University of Munzur, Tunceli, Turkey\\ inanunal@munzur.edu.tr}

\date{ }

\begin{document}

\maketitle

\begin{abstract}
	In this paper the geometry of \NCPM\  is studied under the flatness of conformal, concircular and quasi-conformal curvature tensors. It is proved that a conformal flat \NCPM \ is an Einstein manifold with a negative scalar curvature and has positive sectional curvature. It is also shown that a concircular flat \NCPM\  an Einstein manifold. Finally it is obtained that a quasi-conformally flat \NCPM\  is an Einstein manifold with a positive scalar curvature and, is a space of constant curvature.
\end{abstract}
\section{Introduction}
Contact transformations were defined as a geometric tool to study system of differential equations in 1872 by S.Lie \cite{geiges2001brief}. The subject contact manifold has many applications to other fields of pure mathematics and some applied areas such as mechanics, optics, thermodynamics, or control theory. Also contact manifolds have several application in theoretical physics \cite{kholodenko2013applications}. The Riemannian geometry of contact manifolds give us geometric interpretation about Einstein manifolds which arise from general relativity. Also the curvature concept has center role in Riemannian geometry and curvature properties of manifolds present geometric aspects from algebraic  relations. \par
A conformal transformation is a map which converts a metric to an other with preserving angle between two vector fields.  Conformal curvature tensor is a $ (1,3)- $tensor on a Riemanian manifold which  
is invariant under conformal transformations. This tensor gives important information about the Riemann geometry of the manifold. If it vanishes the manifold is said to be conformally flat, that's mean the manifold is flat under conformal transformations.  A concircular transformation is a special conformal transformation,and preserve the geodesic circle. This type transformations and their applications to differential geometry were studied by Yano \cite{yano1940concircular}.  
Concircular curvature tensor was defined by Yano \cite{yano1940concircular} and it is
invariant under concircular transformations. Also a manifold is called
concircularly flat if this tensor vanishes. Yano and Sawaski \cite{yano1968riemannian} introduced quasi-conformal curvature which 
includes both concircular and conformal curvature as special cases. The Riemannian geometry of contact manifolds is examined with these tensors via flatness and symmetries.  \par
Blair, Ludden and Yano \cite{blair1974geometry} studied on complex manifolds consider the results on Calabi-Eckman manifolds $ S^{2p+1}\times S^{2q+1} $. By consider two Sasakian structure on  $ S^{2p+1}$ and $ S^{2q+1} $ they gave the second fundamental form on Calabi-Eckman manifold, defined Hermitian bicontact manifold and obtained an $ f- $structure on bicontact manifolds. Also normality of bicontact manifolds was given. Bande and Hadjar \cite{bande2005contact} studied on bicontact manifolds under the name contact pairs. Further they considered a special type of f-structure with complementary frames related to a contact pair and called the contact pair structure. Also the normality of contact pair structures were given by same authors \cite{bande2010normal,bande2009contact,bande2013curvature}. \par
The conformal flatness of a \NCPM\ were studied by Bande, Blair and Hadjar \cite{bande2015bochner}. They proved that a conformal flat \NCPM\ is locally isometric to Hopf manifold $ S^{2p+1}(1)\times S^1 $. On the other hand the flatness condition of conformal, concircular and quasi-conformal curvature tensors on contact manifolds has many geometric and physical applications. A conformal flat Sasakian manifold is of constant curvature \cite{de2009complex}. At the same time a normal complex contact metric manifold is not conformal, concircular and quasi-conformal flat \cite{turgut2017conformal}. \par

In this paper we studied on conformal, concircular and quasi-conformal curvature flatness of normal contact pair manifolds. We prove that a conformal flat \NCPM \ is an Einstein manifold with a negative scalar curvature and, has positive sectional curvature. Also we show that a concircular flat \NCPM\  is an Einstein manifold. Finally we prove that a quasi-conformal flat \NCPM\  is an Einstein manifold with a positive scalar curvature and, is a space of constant curvature.

\section{Preliminaries}
In this section a short survey is given for contact manifolds and contact pair structures. For detail about contact manifolds we refer to reader \cite{blair2010riemannian}, and \cite{bande2005contact,bande2010normal,bande2009contact} for contact pairs. 
\subsection{Real and Complex Contact Manifolds}
A real contact manifold is defined by a contact form $ \eta $ which is a volume form on a real $ (2p+1)- $ dimensional differentiable manifold $ M $. The kernel of $ \eta $ defines $ 2p- $dimensional a non-integrable distribution of $ TM $: 
\begin{equation}
\mathcal{D}= \{X:\eta(X)=0, \ X\in\Gamma(TM)\}.
\end{equation}
We also recall $ \mathcal{D} $ contact or horizontal distribution. Let take a vector field $ \xi $ on $ M $ which is dual vector of $ \eta $. Then for $ (1,1)- $tensor field $ \varphi $, $ M $ is called an almost contact metric manifold if it satisfied following conditions: 
\begin{equation}
\varphi^2=-I+\eta\otimes\xi, \ \ \eta(\xi)=1,\ \ , g(\varphi \bullet, \bullet)=-g(\bullet, \varphi \bullet)
\end{equation}
where $ I $ is identity map on $ M $ and $ g $ is a Riemannian metric. Also we recall $ g $ by compatible metric. As the similar to K\"ahler manifold we have a second fundamental form on an almost contact metric manifold $ \Omega(\bullet,\bullet)=d\eta(\bullet,\bullet) $. Also $ d\eta(\bullet,\bullet)=g(\bullet, \varphi \bullet) $ and in this case we recall $ g $ is an associated metric. \\ The geometry of contact manifold is studied in different classes. One of them is Sasakian manifold which has a K\"ahler form on Riemannian cone $ M\times \mathbb{R^+}.  $ A Sasakian manifold has also an almost contact metric structure. The almost contact structure on a Sasakian manifold is normal i.e $ N(\varphi \bullet,\varphi\bullet)+2d\eta(\bullet,\bullet)\xi=0 $ where  $ N(\varphi \bullet,\varphi\bullet)$ is the Nijenhuis tensor field of $ \varphi $.\par
In 1959 Kobayashi \cite{kobayashi1959remarks} defined complex analogue of a real contact manifold. Therefore the concept of complex contact manifold entered to the literature. 1980s Ishihara and Konishi \cite{ishihara1982complex} construct almost contact structure on a complex contact manifold and they defined compatible metric. A complex almost contact metric  manifold  is a complex odd $ (2p+1)- $dimensional complex manifold with $ (J,\varphi,\varphi \circ J,\xi, -J\circ \xi,\eta,\eta \circ J,g) $ structure such that
\begin{eqnarray*}
	&\varphi^2=(\varphi J)^2=-I+\eta\otimes\xi-(\eta \circ J)\otimes (J\circ \xi)\\& \eta(\xi)=1\ \ ,\eta(-J\circ \xi)=0, (\eta \circ J)(-J\circ \xi)=1, \ (\eta \circ J)(\xi)=0\\ &g(\varphi \bullet, \bullet)=-g(\bullet, \varphi \bullet), \ \ g((\varphi \circ J)\bullet,\bullet)=-g(\bullet,(\varphi \circ J)\bullet)
\end{eqnarray*}
where $  g$ is Hermitian metric on $ M $, $ J $ is an almost complex structure. The normality of complex almost contact metric manifolds were given by Ishihara and Konishi and Korkmaz \cite{ishihara1982complex,korkmaz2000}. Normal complex contact metric manifolds were studied by several authors \cite{turgut2017conformal,korkmaz2000,vanli2006boothby}.
\subsection{Metric Contact Pair Manifold}

\begin{dfn}
	Let 	$ M $ be a $ 2p+2q+2 $-dimensional differentiable manifold. A pair of $ (\alpha_1,\alpha_2) $ on $ M $ is said to be a contact pair of type $ (p,q) $ if 
	\begin{itemize}
		\item $ \alpha_1 \wedge (d\alpha_1)^p\wedge\alpha_2 \wedge (d\alpha_2)^q\neq0 $
		\item  $ (d\alpha_1)^{p+1}=0 $ and $ (d\alpha_2)^{q+1}=0 $
	\end{itemize}
	where $ p,q $ are positive integers \cite{bande2005contact}. 
\end{dfn}
For $ 1- $forms $ \alpha_1 $ and $ \alpha_2 $ we have two integrable subbundle of $ TM $; $ \mathcal{D}_{1}=\{X:\alpha_1(X)=0, X\in \Gamma(TM)\} $ and  $ \mathcal{D}_{2}=\{X:\alpha_2(X)=0, X\in \Gamma(TM)\} $. Then we have two characteristic foliations of $ M $, denoted $ \mathcal{F}_1=\mathcal{D}_1 \cap kerd\alpha_1$ and $ \mathcal{F}_2= \mathcal{D}_2 \cap kerd\alpha_2$ respectively. $ \mathcal{F}_1 $ and $ \mathcal{F}_2 $ are $( 2p+1) $ and $ (2q+1)-$dimensional contact manifolds with contact form induced by $ \alpha_1 $ and $ \alpha_2 $. Thus we can define $ (2p+2q)- $dimensional horizontal subbundle $ \mathcal{H} $ 
\begin{equation*}
\mathcal{H}=ker\alpha_1\cap ker\alpha_2.
\end{equation*}
To a contact pair $ (\alpha_1,\alpha_2) $ of type $ (p,q) $ there are associated two commuting vector fields $ Z_1 $ and $ Z_2 $, called Reeb vector fields of the pair, which are uniquely determined by the following equations:
\begin{eqnarray*}
	&\alpha_1(Z_1)=\alpha_2(Z_2)=1, \ \alpha_1(Z_2)=\alpha_2(Z_1)=0\\&
	i_{Z_1}d\alpha_1 =i_{Z_1}d\alpha_2=i_{Z_2}d\alpha_2=0
\end{eqnarray*}
where $ i_X $ is the contraction with the vector field X. In particular, since the Reeb vector fields
commute, they determine a locally free $ \mathbb{R}^2 $-action, called the Reeb action. \\
The tangent bundle of $ (M,(\alpha_1,\alpha_2)) $ can be split into in a different way. For the two subbundle of $ TM $
\begin{eqnarray*}
	T\mathcal{G}_i=ker d\alpha_i\cap ker \alpha_1 \cap ker \alpha_2, \ i=1,2
\end{eqnarray*}
and we can write
\begin{equation*}
T\mathcal{F}_i=T\mathcal{G}_i\oplus \mathbb{R}Z_1.
\end{equation*}
Therefore we get $ TM=T\mathcal{G}_1\oplus T\mathcal{G}_2\oplus \mathbb{R}Z_1\oplus \mathbb{R}Z_2$. The horizontal subbundle can be written as  $ \mathcal{H}=T\mathcal{G}_1\oplus T\mathcal{G}_2 $ and $ \mathcal{V}= \mathbb{R}Z_1\oplus \mathbb{R}Z_2$, we call $ \mathcal{V} $ is  vertical subbundle of $ TM. $ \\
Let $ X $ be an arbitrary vector field on $ M $. We can write $ X=X^{\mathcal{H}}+X^{\mathcal{V}} $, where $ X^{\mathcal{H}}, X^{\mathcal{V}} $ horizontal and vertical component of $ X $ respectively. For $ X^1 \in T\mathcal{F}_1 $ and $ X^2\in T\mathcal{F}_2 $ we have $ X=X^1+X^2 $. Also we can write $ X^1=X^{1^h} +\alpha_2(X^1)Z_2$ and $ X^2=X^{2^h} +\alpha_1(X^2)Z_1$, where $ X^{1^h} $ and $ X^{2^h} $ are horizontal parts of $ X^1 , X^2$ respectively. From all these decomposition of $ X $ finally we get 
\begin{eqnarray*}
	&X=X^{1^h}+X^{2^h}+\alpha_1(X^2)Z_1+\alpha_2(X^1)Z_2\\&\alpha_1(X^{1^h})=\alpha_1(X^{2^h})=0,\ \ \alpha_2(X^{1^h})=\alpha_2(X^{2^h})=0.
\end{eqnarray*}
Since we have two different $ 1- $form by above decomposition we understand the components of $ X\in TM $ in which distributions. 
\begin{dfn}
	An almost contact pair structure on a $ (2p+2q+2)- $dimensional manifold $ M $ is a triple $ \alpha_1, \alpha_2, \phi $, where $ (\alpha_1, \alpha_2) $ is a contact pair and $ \phi $ a $ (1,1) $ tensor field such that:
	\begin{eqnarray}
	\phi^2=-I+\alpha_1\otimes Z_1+\alpha_2\otimes Z_2
	,\ \ \phi Z_1=\phi Z_2=0. 
	\end{eqnarray}
	The rank of $ \phi $ is $ (2p+2q) $ and $ \alpha_i(\phi)=0 $ for $ i=1,2 $. 
\end{dfn}
The endomorphism $\phi $ is said to be decomposable if $ T\mathcal{F}_i $ is invariant under $ \phi $. If $ \phi $ is decomposable then  $ (\alpha_i, Z_i, \phi) $ induce an almost contact structure on $ \mathcal{F} _i$ for $ i=1,2 $ \cite{bande2005contact}. Unless otherwise stated we assume that $ \phi $ is decomposable. 
\begin{dfn}
	Let $ (\alpha_1,\alpha_2,Z_1,Z_2,\phi) $ be an almost contact pair structure on a manifold $ M $. A Riemannian metric $ g $ is called 
	\begin{enumerate}
		\item compatible if $ g(\phi X_1,\phi X_2)=g(X_1,X_2)-\alpha_1(X_1)\alpha_1(X_2)-\alpha_2(X_1)\alpha_2(X_2) $ for all $ X_1,X_2 \in TM $,
		\item  associated if $ g(X_1,\phi X_2)=(d\alpha_1+d\alpha_2)(X_1,X_2) $ and $ g(X_1,Z_i)=\alpha_i(X_1) ,$ for $ i=1,2 $ and for all $ X_1,X_2 \in \Gamma(TM) $
	\end{enumerate}
\end{dfn}
$ 4- $tuple $ (\alpha_1,\alpha_2,\phi, g) $ is called a metric almost contact pair on a manifold $ M $ and g is an associated metric with respect to contact pair structure $ (\alpha_1,\alpha_2,\phi) $. We recall $ (M,\phi,Z_1,Z_2,\alpha_1,\alpha_2,g) $ is a metric contact pair manifold.\\
We have following properties for a metric almost contact pair manifold $ M $ \cite{bande2005contact}: 
\begin{eqnarray}
&g(Z_i,X)=\alpha_i(X), \ \  g(Z_i,Z_j)=\delta_{ij} \\& \nabla_{Z_i}Z_j=0, \ \nabla_{Z_i}\phi=0
\end{eqnarray}
and for every $ X $ tangent to $ \mathcal{F}_i \,\ i=1,2$ we have
\begin{equation*}
\nabla_ X Z_1=-\phi_1  X , \ \ \nabla_ X Z_2=-\phi_2  X 
\end{equation*}
where $ \phi=\phi_1+\phi_2 $. Normality of metric contact pair manifold is given by Bande and Hadjar \cite{bande2010normal}. Let  $ (M,\phi,Z_1,Z_2,\alpha_1,\alpha_2,g) $ be metric contact pair manifold then we have two almost complex structure: 
\begin{equation*}
J=\phi-\alpha_2\otimes Z_1+\alpha_1\otimes Z_2,\ \ T=\phi+\alpha_2\otimes Z_1-\alpha_1\otimes Z_2
\end{equation*}
\begin{dfn}
	A metric contact pair manifold is said to be normal if $ J $ and $ T $ are integrable \cite{bande2010normal}. 
\end{dfn}
\begin{teo}
	Let $ (M,\phi,Z_1,Z_2,\alpha_1,\alpha_2,g) $ be a normal metric contact pair manifold then we have
	\begin{equation}\label{Phicovariantderivation}
	g((\nabla_{X_1}\phi)X_2,X_3)=\sum_{i=1}^{2}(d\alpha_i(\phi X_2,X_1)\alpha_i(X_3)-d\alpha_i(\phi X_3,X_1)\alpha_i(X_2))
	\end{equation}
	where $ X_1,X_2,X_3 $ are arbitrary vector fields on $ M $ \cite{bande2009contact}.
\end{teo}
We can consider a natural question: could any metric contact pair structure be considered locally the product of two contact metric manifold? An example of metric contact pair were given in \cite{bande2011characteristic}, which is not locally product of two contact metric manifold. So metric contact pair structure has some different properties from contact metric manifolds and their results will have useful interpretation for the geometry of contact and complex manifolds.
\\
On a normal metric contact pair manifold we have $ \nabla_XZ=-\phi X $ for $ X\in \Gamma(TM) $ and $ Z=Z_1+Z_2 $. 

\subsection{Curvature Properties of Normal Metric Contact Pair Manifolds}
We  use the following statements for the Riemann curvature;
\begin{eqnarray*}
	&R(X_1,X_2)X_3=\nabla_{X_1}\nabla_{X_2}X_3-\nabla_{X_2}\nabla_{X_1}X_3-\nabla_{[X_1,X_2]}X_3,\\& \mathcal{R}(X_1,X_2,X_3,X_4)=g(R(X_1,X_2)X_3,X_4).
\end{eqnarray*}
for all $ X_1,X_2,X_3,X_4 \in \Gamma(TM) $. Also the Ricci operator is defined by
\begin{equation*}
QX=\sum_{i=1}^{dim(M)}R(X,E_i)E_i,
\end{equation*} 
the Ricci curvature and scalar curvature is given by 
\begin{eqnarray*}
	Ric(X_1,X_2)=g(QX_1,X_2),\\scal=\sum_{i=1}^{dim(M)}Ric(E_i,E_i).
\end{eqnarray*}

In \cite{bande2013curvature} the curvature of contact pairs were examined. Let $(M,\phi,Z_1,Z_2,\alpha_1,\alpha_2,g)$ be a normal metric contact pair manifold. Then for  $ X_1,X_2, X_3 \in \Gamma(\mathcal{H}) $ and $ Z=Z_1+Z_2 $ Reeb vector field we have :
\begin{eqnarray}
R(X_1,Z)X_2&=&-d{\alpha_1}(\phi X_1,X_2)Z_1-d{\alpha_2}(\phi X_1,X_2)Z_2 \\&=&-[(d\alpha_1+d\alpha_2)(\phi X_1,X_2)]Z, \notag
\end{eqnarray}
\begin{eqnarray}\label{R(X,Y,Z,V)}
\mathcal{R}(X_1,X_2,Z,X_3)&=d\alpha_1(\phi X_3,X_1)\alpha_1(X_2)+d\alpha_2(\phi X_3, X_1)\alpha_2(X_2)\\&-d\alpha_1(\phi X_3,X_1)\alpha_1(X_2)-d\alpha_2(\phi X_3,X_2)\alpha_2(X_1) \notag,
\end{eqnarray}
\begin{equation}\label{R(X,Z)Z}
R(X_1,Z)Z=-\phi^2X_1.
\end{equation}
Let take an orthonormal basis of $ M $ 
\begin{equation*}
\{E_1,E_2,...,E_p,\phi E_1,\phi E_2,...,\phi E_p,E_{p+1},E_{p+2},...,E_{p+q}, \phi E_{p+1},\phi E_{p+2},...,\phi E_{p+q} , Z_1,Z_2\}
\end{equation*}
then for all $ X_1\in \Gamma(TM) $ we get the Ricci curvature of $ M $ as 
\begin{equation}
Ric(X_1,Z)=\sum_{i=1}^{2p+2q}{d\alpha_1(\phi E_i, E_i)\alpha_1(X)+d\alpha_2(\phi E_i, E_i)\alpha_2(X)}.
\end{equation}
 So, we obtain the following result:
\begin{eqnarray}
&Ric(X,Z)=0,\ \ \text{for } X\in \Gamma(\mathcal{H}) \label{Rİc(X,Z)},\\
&Ric(Z,Z)=2p+2q.\label{Ric(Z,Z)}\\
&Ric(Z_1,Z_1)=2p , \ Ric(Z_2,Z_2)=2q , \ Ric(Z_1,Z_2)=0.
\end{eqnarray}
Conformal $ \mathcal{C} $, concircular $ \mathcal{W} $ and quasi-conformal curvature tensor $ \widetilde{\mathcal{C}}$ of a $(2p+2q+2)$-dimensional normal contact metric pair manifold are given by followings, respectively: 
\begin{eqnarray}
\mathcal{C}( X_1,X_2) X_3 &=&R( X_1,X_2) X_3+\frac{scal}{(
	2p+2q+1) (2p+2q)}( g(X_2,X_3)X_1-g(X_1,X_3)X_2)   \label{WEX_2L} \\
&&+\frac{1}{2p+2q}( g( X_1,X_3) QX_2-g( X_2,X_3) QX_1+Ric(
X_1,X_3) X_2-Ric( X_2,X_3) X_1) , \notag\\\notag \\
\mathcal{W}( X_1,X_2) X_3&=&R( X_1,X_2) X_3-\frac{scal}{(
	2p+2q+2) ( 2p+2q+1) }[ g(X_2,X_3)X_1-g(X,X_3)X_2] ,
\label{concircular} \\\notag\\
\widetilde{\mathcal{C}}( X_1,X_2) X_3 &=&aR( X_1,X_2) X_3+b[
Ric( X_2,X_3) X_1-Ric( X_1,X_3) X_2
\label{quasi conf.tensor} \\
&&+g( X_2,X_3) QX_1 -g( X_1,X_3) QX_2] \notag\\
&&-\frac{scal}{2p+2q+2}[ \frac{a}{2p+2q+1}%
+2b] [ g(X_2,X_3)X_1  \notag -g(X_1,X_3)X_2]  \notag
\end{eqnarray}
where  $ X_1,X_2,X_3 \in \Gamma(TM) $,  $a$ and $b$ are constants. 

\section{Hermitian Contact Pair Manifold}
As known, the product of two contact metric manifolds is a contact pair metric manifold. In this section we give an almost contact pair structure on a Hermitian manifold.\\
Let $ (M,g,J) $ be $ (2p+2q+2)- $dimensional Hermitian manifold and $ (\varphi_1, \eta_1, \xi_1) $ and $ (\varphi_2, \eta_2, \xi_2) $ be two almost contact structures on $ M $ with following properties.
\begin{eqnarray*}
	&g(\varphi_i X_1,X_2)=-g(X_1,\varphi_i X_2), \ \ \text{for  } i=1,2\\& J\xi_1=\xi_2,\ \ J\xi_2=-\xi_1 \\ & {\varphi_i}^2X_1=-X_1+\eta_1(X_1)\xi_1+\eta_2(X_1)\xi_2
	\\& \varphi_1(JX_1)=-J\varphi_1X_1=\varphi_2X_1\\ &\varphi_2(JX_1)=-J\varphi_2X_1=-\varphi_1X_1\\& \varphi_2(\varphi_1X_1)=-\varphi_1(\varphi_2X_1)=JX_1+\eta_1(X_1)\xi_2-\eta_2(X_1)\xi_1
\end{eqnarray*}
where $ X_1,X_2  $ are two arbitrary vector fields on $ M $ \cite{beldjilali2016structures}. 

Let take $ \phi=\varphi_1 \circ \varphi_2 $. Then $ \phi $ is a $ (1,1) $ tensor field on $ M $. By direct computation we get 
\begin{equation*}
\phi^2X_1=-X_1+\eta_1(X_1)\xi_1+\eta_2(X_1)\xi_2.
\end{equation*}
Thus we obtain an almost contact pair structure on $ M $ with the contact pair $ \eta_1, \eta_2 $ and we state:
\begin{cor}
	Let $ (M^{2p+2q+2}, J,g) $ be an almost Hermitian manifold and $ (\varphi_i, \eta_i,\xi_i)_{i=0}^{2} $ be two almost contact structure on $ M $ with properties are given above. Then $ (\eta_1, \eta_2, \phi) $ is an almost contact pair structure on $ M $ such that 
	\begin{eqnarray*}
		&\phi^2X_1=-X_1+\eta_1(X_1)\xi_1+\eta_2(X_1)\xi_2\\ & \eta_i(\xi_j)=\delta_{ij}, \ 1\leq i,j,\leq 2\\&
		\phi(\xi_i)=0
	\end{eqnarray*}
	for all $ X_1 \in \Gamma(TM)$. 
\end{cor}
Also for $ X_1,X_2 \in \Gamma (TM) $ we have 
\begin{equation*}
g(\phi X_1,X_2)=-g(X_1,\phi X_2)
\end{equation*}
and 
\begin{equation*}
g(\phi X_1,\phi X_2)=g(X_1,X_2)-\eta_1(X_1)\eta_1(X_2)-\eta_2(X_1)\eta_2(X_2).
\end{equation*}
Thus we obtain compatible metric with contact pair structure. \par 
These results show that a contact pair structure on an almost Hermitian manifold could be obtained from two almost contact structure on this manifold. Since contact pair manifolds have some significant properties, some future works could be done for Hermitian and contact structure. Also if the manifold is complex the structure $ (\varphi_1,\varphi_2,J,\xi_1,\xi_2,\eta_1,\eta_2,g) $ be a complex almost contact metric manifold. This type of manifolds was studied by several authors \cite{turgut2017conformal, kobayashi1959remarks,ishihara1982complex,korkmaz2000,vanli2006boothby}

\section{Flatness Conditions on Normal Contact Pair Manifold}
In this section we give some results on the flatness of conformal, concircular and quasi-conformal curvature tensors. 
\begin{teo}
	A conformal flat \NCPM \ is an Einstein manifold with a negative scalar curvature and has positive sectional curvature. 
\end{teo}

\begin{prf}
	Let $ (M,\phi,\alpha_1,\alpha_2) $ be \NCPM. Suppose that $ M $ is conformal flat. Then we have 
	\begin{eqnarray}\label{conformalproof1}
	\mathcal{R}(X_1,X_2,X_3,X_4)&=-A[g(X_2,X_3)g(X_1,X_4)-g(X_1,X_3)g(X_2,X_4) ]\\& -B( g( X_1,X_3) Ric(X_2,X_4)-g( X_2,X_3) Ric(X_1,X_4)\notag\\&+Ric(
	X_1,X_3) g(X_2,X_4)-Ric( X_2,X_3)Ric( X_1,X_4)).\notag
	\end{eqnarray} 
	where $ A=\frac{scal}{(2p+2q+1)(2p+2q)} $ and $ B= \frac{1}{2p+2q}$. 
	Taking $ X_2=X_3=Z $ and $ X_1,X_4 \in \Gamma(\mathcal{H})$ in (\ref{conformalproof1}) , since $ g(Z,Z)=2 $ and from (\ref{R(X,Z)Z}), (\ref{Rİc(X,Z)})  we obtain  
	\begin{eqnarray}
	\mathcal{R}(X_1,Z,Z,X_4)=-2Ag(X_1,X_4)-2BRic(X_1,X_4)
	\end{eqnarray}
	Also from (\ref{R(X,Y,Z,V)}) we get 
	\begin{equation*}
	(d\alpha_1+\alpha_2)(\phi X_4,X_1)=-2Ag(X_1,X_4)-2BRic(X_1,X_4)
	\end{equation*}
	and thus we obtain
	\begin{equation}\label{Einstein1}
	Ric(X_1,X_4)=-\frac{2A+1}{2B}g(X_1,X_4).
	\end{equation}
	So, the manifold is Einstein. On the other hand by direct computation from (\ref{Einstein1}) the scalar curvature is
	\begin{equation}\label{scalar1}
	scal=-\frac{(2p+2q+1)(2p+2q+2)}{2p+2q+3}.
	\end{equation}
	This shows the scalar curvature is negative. Let choose $ X_1=X_4 $, $ X_2=X_3 $ unit and orthogonal vector fields in (\ref{conformalproof1}). Then  the sectional curvature is obtained by 
	\begin{equation*}
	k(X_1,X_2)=-A=-\frac{(2p+2q+1)^2(2p+2q+2)(2p+2q)}{2p+2q+3}.
	\end{equation*}
	Thus, the proof is completed.
\end{prf}

An Einstein manifold is also Einstein under concircular transformation. Yano proved that a concircular flat Riemann manifold is Einstein\cite{yano1940concircular}. By similar  way we can easily obtain following result. 
\begin{teo}
	A concircular flat normal contact pair manifold is Einstein. 
\end{teo}
Our finally result is about quasi-conformal flatness of \NCPM. 
\begin{teo}
	A quasi-conformally flat \NCPM; 
	\begin{enumerate}
		\item is an Einstein manifold with a positive scalar curvature
		\item is a space of constant curvature.
	\end{enumerate}
\end{teo}

\begin{prf}
	Let $ M $ be a quasi-conformally flat \NCPM. Then for $ X_1,X_2,X_3, X_4 \in \Gamma(TM) $ we have 
	\begin{eqnarray}
	0=&&a\mathcal{R}( X_1,X_2,X_3,X_4) +b[
	Ric( X_2,X_3)g( X_1,X_4)-Ric( X_1,X_3) g(X_2,X_4)
	\label{quasi conf.tensor.in.proof} \\
	&&+g( X_2,X_3) Ric(X_1,X_4) -g( X_1,X_3) g(X_2,X_4)]\notag \\
	&& -\frac{scal}{2p+2q+2}[ \frac{a}{2p+2q+1}%
	+2b] [ g(X_2,X_3)g(X_1,X_4)  -g(X_1,X_3)g(X_2,X_4)]  \notag
	\end{eqnarray}
	Let write  $\frac{scal}{2p+2q+2}[ \frac{a}{4m+1}+2b]=K $ for brevity. In ( \ref{quasi conf.tensor.in.proof}), by taking $ X_1=X_4=E_i $ and getting sum from $ i=1 $ to $ i=2p+2q+2 $ we obtain 
	\begin{eqnarray*}
		0=(a+b(2p+2q))Ric(X_2,X_3)+(b)scal-K(2p+2q+1))g(X_2,X_3)
	\end{eqnarray*}
	and therefore we get
	\begin{equation*}
	0=[a+b(2p+2q)][Ric(X_2,X_3)-\frac{scal}{2p+2q+2}g(X_2,X_3)]. 
	\end{equation*}
	Assume that $ a+b(2p+2q)\neq 0 $. Then 
	\begin{equation}\label{quasiconfprof2}
	Ric(X_2,X_3)=\frac{scal}{2p+2q+2}g(X_2,X_3).
	\end{equation}
	By taking $ X_2=X_3=Z $ in (\ref{quasiconfprof2}) and from (\ref{Ric(Z,Z)}), we get positive scalar curvature is given by
	\begin{equation*}
	scal=4(p+q)(p+q+1).
	\end{equation*}
	So, the Ricci curvature has the following form:
	\begin{equation}\label{quasiconfprof3}
	Ric(X_2,X_3)=2(p+q)g(X_2,X_3).
	\end{equation}
	This shows manifold is Einstein.\\
	On the other hand consider (\ref{quasiconfprof3}) in (\ref{quasi conf.tensor.in.proof}) we get 
	\begin{eqnarray*}
		0=&&a[\mathcal{R}(X_1,X_2,X_3,X_4)-\frac{p+q}{2p+2q+1}[g(X_2,X_3)g(X_1,X_4)\\&&-g(X_1,X_3)g(X_2,X_4)]
	\end{eqnarray*}
	If $ a\neq 0 $ we get 
	\begin{eqnarray*}
		\mathcal{R}(X_1,X_2,X_3,X_4)=\frac{p+q}{2p+2q+1}[g(X_2,X_3)g(X_1,X_4)-g(X_1,X_3)g(X_2,X_4)].
	\end{eqnarray*}
	This shows us that the manifold is a space of constant curvature.
	
\end{prf}




\vspace{6pt}


\begin{thebibliography}{999}
	\bibitem{geiges2001brief}Geiges, H. (2001). A brief history of contact geometry and topology. Expositiones Mathematicae, 19(1), 25-53.
	\bibitem{kholodenko2013applications} Kholodenko, A. L. (2013). Applications of contact geometry and topology in physics. World Scientific.
	\bibitem{yano1940concircular}
	Yano, K. (1940). Concircular geometry I. Concircular transformations. Proceedings of the Imperial Academy, 16(6), 195-200.
	\bibitem{yano1968riemannian} Yano K. and Sawaski S.(1968). Riemannian manifolds admitting a
	conformal transformation group. J. Diff. Geo. 2 \ , 161-184
	\bibitem {blair1974geometry}Blair, D. E., Ludden, G. D.,  Yano, K. (1974). Geometry of complex manifolds similar to the Calabi-Eckmann manifolds. Journal of Differential Geometry, 9(2), 263-274.
	\bibitem {bande2005contact} Bande, G.,  Hadjar, A. (2005). Contact pairs. Tohoku Mathematical Journal, Second Series, 57(2), 247-260.
	\bibitem{bande2010normal} Bande, G., Hadjar, A. (2010). On normal contact pairs. International Journal of Mathematics, 21(06), 737-754.
	\bibitem{bande2009contact}Bande, G.,  Hadjar, A. (2009). Contact pair structures and associated metrics. In Differential Geometry (pp. 266-275).	\bibitem{bande2013curvature}Bande, G., Blair, D. E.,  Hadjar, A. (2013). On the curvature of metric contact pairs. Mediterranean journal of mathematics, 10(2), 989-1009.
	
	\bibitem{bande2015bochner}Bande, G., Blair, D. E., Hadjar, A. (2015). Bochner and conformal flatness of normal metric contact pairs. Annals of Global Analysis and Geometry, 48(1), 47-56.
	\bibitem{de2009complex} De, U. C.,  Shaikh, A. A. (2009). Complex manifolds and contact manifolds. Narosa Publishing House.
	\bibitem{turgut2017conformal} Turgut Vanli, A., Unal, I. (2017). Conformal, concircular, quasi-conformal and conharmonic flatness on normal complex contact metric manifolds. International Journal of Geometric Methods in Modern Physics, 14(05), 1750067.
	\bibitem{blair2010riemannian} Blair, D. E. (2010). Riemannian geometry of contact and symplectic manifolds. Springer Science Business Media.
	\bibitem{kobayashi1959remarks} Kobayashi, S., Remarks on complex contact manifolds, Proc.
	Amer. Math. Soc. v.10, 164--167 (1959).
	\bibitem{ishihara1982complex} Ishihara, S. and Konishi, M. (1982). Complex almost contact
	structures in a complex contact manifold, Kodai Math. J. v.5, 30--37 .
	\bibitem{korkmaz2000} Korkmaz, B. (2000). Normality of complex contact manifolds, Rocky Mountain J. Math. v.30, 1343--1380 
	\bibitem{vanli2006boothby} Turgut Vanli, A. and Blair, D. E., The Boothby-Wang Fibration of the Iwasawa Manifold as a Critical Point of the Energy, Monatsh. Math. v.147, 75--84 (2006).\bibitem{bande2011characteristic} Bande, G., Hadjar, A. (2011). On the characteristic foliations of metric contact pairs. Harmonic Maps and Differential Geometry. Contemp. Math, 542, 255-259.
	
	\bibitem{beldjilali2016structures}Beldjilali, G., Belkhelfa, M. (2016). Structures on the Product of Two Almost Hermitian Almost Contact Manifolds. International Electronic Journal of Geometry, 9(2), 80-86.

	
	
	
\end{thebibliography}
\end{document}